\documentclass[reqno,10pt]{amsart}
\usepackage{amscd,amssymb}
\usepackage{color}
\usepackage[colorlinks=true,linkcolor=blue,citecolor=blue]{hyperref}
\numberwithin{equation}{section}
\theoremstyle{plain}

\newcommand\eps{\varepsilon}
\newcommand\zet{\zeta}
\newcommand\tet{\theta}

\newcommand\lam{\lambda}                
         \newcommand\Sig{\Sigma}

\newcommand\ome{\omega}         

\newcommand\calF{{\mathcal{F}}}

\newcommand\calM{{\mathcal{M}}}

\newcommand\calP{{\mathcal{P}}}
\newcommand\calQ{{\mathcal{Q}}}

\newcommand\RR{\mathbb{R}}

\newcommand\PP{\mathbb{P}}

\newcommand\ZZ{\mathbb{Z}}

\newcommand\CC{\mathbb{C}}

 \newcommand\grt{{\mathfrak{t}}}

\newcommand\nek{,\ldots,}
\newcommand\sdp{\times \hskip -0.3em {\raise 0.3ex
\hbox{$\scriptscriptstyle |$}}} 

\newcommand\const{\operatorname{const}}

\newcommand\diag{\operatorname{diag}}

\newcommand\End{\operatorname{End\,}}

\newcommand\Hom{\operatorname {Hom}}

\newcommand\ind{\operatorname{ind}}

\newcommand\Tr{\operatorname{Tr}}

\newcommand\oz{{\overline{z}}}

\newcommand\tilT{{\widetilde{T}}}

\renewcommand{\>}{\rangle}
\newcommand{\<}{\langle}

\newcommand{\Cech}{\v{C}ech\ }

\theoremstyle{plain}
\newtheorem{Thm}[subsection]{Theorem}
\newtheorem{Cor}[subsection]{Corollary}
\newtheorem{Lem}[subsection]{Lemma}
\newtheorem{Prop}[subsection]{Proposition}
\newtheorem{Conjec}[subsection]{Conjecture}

\theoremstyle{definition}
\newtheorem{Def}[subsection]{Definition}

\theoremstyle{remark}

\newtheorem{Rem}[subsection]{Remark}

\newif\ifShowLabels
\ShowLabelstrue
\newdimen\theight
\def\TeXref#1{%
        \leavevmode\vadjust{\setbox0=\hbox{{\tt
                \quad\quad  {\small \textrm #1}}}%
        \theight=\ht0
        \advance\theight by \lineskip
        \kern -\theight \vbox to
        \theight{\rightline{\rlap{\box0}}%
        \vss}%
        }}%

\ShowLabelsfalse

\renewcommand{\sec}[2]{\section{#2}\label{S:#1}%
        \ifShowLabels \TeXref{{S:#1}} \fi}
\newcommand{\ssec}[2]{\subsection{#2}\label{SS:#1}%
        \ifShowLabels \TeXref{{SS:#1}} \fi}

\newcommand{\refs}[1]{Section ~\ref{S:#1}}
\newcommand{\refss}[1]{Subsection ~\ref{SS:#1}}
\newcommand{\reft}[1]{Theorem ~\ref{T:#1}}
\newcommand{\refl}[1]{Lemma ~\ref{L:#1}}
\newcommand{\refp}[1]{Proposition ~\ref{P:#1}}
\newcommand{\refc}[1]{Corollary ~\ref{C:#1}}

\newcommand{\refr}[1]{Remark ~\ref{R:#1}}
\newcommand{\refe}[1]{\eqref{E:#1}}

\newenvironment{thm}[1]%
        { \begin{Thm} \label{T:#1}  \ifShowLabels \TeXref{T:#1} \fi }%
        { \end{Thm} }

\renewcommand{\th}[1]{\begin{thm}{#1} \sl }
\renewcommand{\eth}{\end{thm} }

\newenvironment{lemma}[1]%
        { \begin{Lem} \label{L:#1}  \ifShowLabels \TeXref{L:#1} \fi }%
        { \end{Lem} }
\newcommand{\lem}[1]{\begin{lemma}{#1} \sl}
\newcommand{\elem}{\end{lemma}}

\newenvironment{propos}[1]%
        { \begin{Prop} \label{P:#1}  \ifShowLabels \TeXref{P:#1} \fi }%
        { \end{Prop} }
\newcommand{\prop}[1]{\begin{propos}{#1}\sl }
\newcommand{\eprop}{\end{propos}}

\newenvironment{corol}[1]%
        { \begin{Cor} \label{C:#1}  \ifShowLabels \TeXref{C:#1} \fi }%
        { \end{Cor} }
\newcommand{\cor}[1]{\begin{corol}{#1} \sl }
\newcommand{\ecor}{\end{corol}}

\newenvironment{conjec}[1]%
        { \begin{Conjec} \label{Conj:#1}  \ifShowLabels \TeXref{C:#1} \fi }%
        { \end{Conjec} }
\newcommand{\conj}[1]{\begin{conjec}{#1} \sl }
\newcommand{\econj}{\end{conjec}}

\newenvironment{defeni}[1]%
        { \begin{Def} \label{D:#1}  \ifShowLabels \TeXref{D:#1} \fi }%
        { \end{Def} }
\newcommand{\defe}[1]{\begin{defeni}{#1} \sl }
\newcommand{\edefe}{\end{defeni}}

\newenvironment{remark}[1]%
        { \begin{Rem} \label{R:#1}  \ifShowLabels \TeXref{#1} \fi }%
        { \end{Rem} }
\newcommand{\rem}[1]{\begin{remark}{#1}}
\newcommand{\erem}{\end{remark}}

\newcommand{\Label}[1]%
        {\label{#1}%
            \ifShowLabels \TeXref{#1} \fi }%

\newcommand{\eq}[1]%
         {\newline \ifShowLabels \TeXref{{E:#1}} \fi
                \begin{equation} \label{E:#1}}
\newcommand{\eeq}{\end{equation}}

\newcommand{\prf}{ \begin{proof} }
\newcommand{\eprf}{ \end{proof} }

\newcommand{\ka}{K\"ahler }
\newcommand\tiltet{\widetilde{\tet}}
\newcommand\F{\calF}
\newcommand\fo{\F{|_Z}\otimes o(Z)}
\begin{document}

\title{Cohomology of a Hamiltonian $T$-space with involution}
\author{Semyon Alesker$^\dag$} 
\thanks{${}^\dag$Supported in part by the ISF grant 1447/12 .}
\address{School of Mathematical Sciences, Tel Aviv University
Ramat Aviv, 69978 Tel Aviv, Israel}

\email{semyon@post.tau.ac.il}
\thanks{${}^\ddag$Supported in part by the NSF grant DMS-1005888.}
\author{Maxim Braverman$^\ddag$}

\address{Department of Mathematics,
        Northeastern University,
        Boston, MA 02115,
        USA
         }

\email{maximbraverman@neu.edu}  

\begin{abstract}
Let $M$ be a compact symplectic manifold on which a compact torus $T$ acts Hamiltonialy with a moment map $\mu$. Suppose there exists a symplectic involution $\tet:M\to M$, such that $\mu\circ\tet=-\mu$. Assuming that 0 is a regular value of $\mu$, we calculate the character of the action of $\tet$ on the cohomology of $M$ in terms of the trace of the action of $\tet$ on the symplectic reduction $\mu^{-1}(0)/T$ of $M$.  This result generalizes a theorem of R.~Stanley, who considered the case when $M$ was a toric variety and $\dim T=\frac12\dim_\RR{}M$.
\end{abstract}

\maketitle

\sec{introd}{Introduction}

In \cite{Stanley87}, R.~Stanley proved a lower bound for the number of faces of a centrally symmetric simple polytope. The main ingredient of his proof is the following result: Suppose $M$ is a toric variety corresponding to a simple centrally symmetric polytope $P$ of dimension $n$. Then $M$ is acted upon by an $n$-dimensional torus $T$. The symmetry on $P$ induces an involution $\tet:M\to M$, such that
\eq{compat0}
      \tet\circ t \ = \ t^{-1}\circ \tet, \qquad \text{for any} \quad
                         t\in T.
\end{equation}
Then $\tet$ acts on the cohomology $H^{2i}(M,\CC)$ of $M$. The result of Stanley states that the trace of this action is equal to $\binom{n}{i}$.

In this paper, we extend the result of Stanley to an arbitrary compact
symplectic manifold $M$, which possesses a Hamiltonian action of a torus $T$ of
arbitrary dimension and an involution $\tet$, which preserves the symplectic
form and satisfies \refe{compat0}. In this situation, the moment map for the
action of $T$ can be chosen so that $\mu\circ\tet=-\mu$. Assuming that $0$ is a
regular value of $\mu$, we compute (cf. \reft{main}) the character of the action
of $\tet$ on the cohomology $H^*(M,\CC)$ in terms of the action of $\tet$ on the
symplectic reduction $M_0=\mu^{-1}(0)/T$ of $M$. If
$\dim{}T=\frac12\dim{}M$, so that $M$ is a toric variety, then $M_0$ is a point
and our formula reduces to the theorem of Stanley. In particular we obtain  a new more geometric proof of the Stanley's result.

The proof of our main theorem is based on a study of the action of
$\tet$ on the equivariant cohomology of $M$. We compute
the {\em graded character} $\chi_\tet(t)$ of this action (cf.
\refe{chi}) in two different ways.

First, we uses the well known formula $H^*_T(M)=H^*(M)\otimes H_T(pt)$
to express $\chi_\tet(t)$ in terms of the character of the action of
$\tet$ on $H^*(M)$, cf. \refp{HT=HxA}.

Our second computation (cf. \refp{chitet}) expresses $\chi_\tet(t)$ in terms of the character of the action of $\tet$ on the cohomology of the symplectic
reduction $H^*(M_0)$. To this end we use equivariant Morse-type inequalities for the square
$f=\<\mu,\mu\>$ of the moment map. This function has degenerate singularities. In general, the set $C$ of singular points of $f$ is not even a manifold. However, Kirwan \cite{Kirwan84} constructed Morse-type inequalities for this type of functions. We will refer to functions satisfying the Kirwan theorem as {\em Kirwan-Morse} functions. Thus $f=\<\mu,\mu\>$ is a Kirwan-Morse function on $M$ equivariant with respect to the action of the semi-direct product
$\ZZ_2\ltimes{}T=\tilT$ on $M$ (here the action of $\ZZ_2$ is generated by $\tet$). In \refl{tetBT} we show that the $\tilT$-equivariant Morse inequalities for $f$ {\em with local coefficients} are, in fact, equalities and, hence, may be used to calculate the equivariant cohomology of $M$. This generalizes a result of Kirwan, \cite[Th.~5.4]{Kirwan84}. Comparing the above Morse {\em equalities} with different local coefficients, we calculate $\chi_\tet(t)$ via the character of the action of $\tet$ on $M_0$.

Comparison of the above two expressions for $\chi_\tet(t)$ leads to a proof of our main theorem.

\subsection*{Contents}
The paper is organized as follows:

In \refs{main} we formulate and prove our
main result -- \reft{main}. The proof is based on two statements
(Propositions~\ref{P:HT=HxA} and \ref{P:chitet}), which are proven in
the later sections.

In \refs{Ex} we present some examples and applications of
\reft{main}. In particular, we show that it implies the result of
Stanley \cite{Stanley87}. We also discuss applications of our theorem
to flag varieties.

In \refs{HT=HxA} we study the action of the involution $\tet$
on the equivariant cohomology of $M$ and prove \refp{HT=HxA} (the
expression of $\chi_\tet(t)$ in terms of the action of $\tet$ on $H^*(M)$).

Finally, in \refs{Morse}, we use the equivariant Morse inequalities to
prove \refp{chitet} (the expression of $\chi_\tet(t)$ in terms of the
action of $\tet$ on $H^*(M_0)$).

\sec{main}{Main theorem}

\ssec{setting}{}
Let $(M,\ome)$ be a compact $2n$-dimensional symplectic manifold endowed with a
Hamiltonian action of a compact $k$-dimensional torus $T$. Let
$\grt^*\simeq\RR^k$ denote the dual space to the Lie algebra of $T$
and let $\mu:M\to \grt^*$ be the moment map for the action of $T$ on
$M$.

Let $\tet:M\to M$ be an involution of $M$ such that  $\tet^*\ome=\ome$
and
\eq{compat}
        \tet (t\cdot m) \ = \ (t^{-1})\cdot m,
                \quad \text{for any} \quad t\in T,  \ m\in M.
\end{equation}
We  will normalize the moment map (which is defined up to an additive
constant) so that
\eq{normal}
     \mu\circ\tet=-\mu.
\end{equation}
In this situation $\tet$ acts naturally on the cohomology $H^i(M)=H^i(M;\CC)$ of
$M$ with complex coefficients. Let $H^i(M)^+$ (resp. $H^i(M)^-$) denote the
subspace of $H^i(M)$ fixed by $\tet$ (resp. the subspace of $H^i(M)$ on which
$\tet$ acts as a multiplication by $-1$). Set
\[
      h^{i,\pm}=\dim_\CC{}H^i(M)^\pm.
\]

Suppose now that zero is a regular value for the moment map $\mu$. Then
$\mu^{-1}(0)$ is a smooth manifold which is invariant under the actions of $T$
and $\tet$. Moreover, the action of $T$ on $\mu^{-1}(0)$ is {\em locally free},
i.e. each point of $\mu^{-1}(0)$ has at most finite stabilizer in $T$. Hence,
the symplectic reduction $M_0:=\mu^{-1}(0)/T$ is an {\em orbifold}. Let
$H^i(M_0)$ denote the cohomology of $M_0$ with complex coefficients.

The involution $\tet$ preserves $\mu^{-1}(0)$ and, hence, acts on $M_0$ and
$H^i(M_0)$. Let $h^{i,+}_0$ (resp. $h^{i,-}_0$) denote the dimension of the
subspace of $\tet$ invariant vectors in $H^i(M_0)$ (resp. the dimension of the
subspace of the vectors on which $\tet$ acts by multiplication by $-1$).

Our principal result is the following
\th{main}
In the situation described above
\eq{main}
        h^{i,+} \ - \ h^{i,-} \ = \
          \sum_{j=0}^{\min\, (k,[i/2])}
                \binom{k}{j}\, (h^{i-2j,+}_0 - h^{i-2j,-}_0),
\end{equation}
for any $i=0,1\nek 2n$.
\eth
\rem{main}
a. \ The conditions of the theorem imply that $\tet$ acts freely on
the set $C$ of fixed points of any subtorus $T'$ of $T$. Moreover,
$\tet$ acts freely on the set of connected components of $C$. It
follows from the fact that,
if $\mu'$ is the moment map for the $T'$-action, then
$\mu'\circ\tet=-\mu'$ and the restriction of $\mu'$ on any
connected component is a non-zero constant.

b. \ The theorem remains true if $M$ is a symplectic orbifold, rather than a
smooth manifold. The proof is just a bit more complicated than the one we
present is this paper.

c. \ Moreover, the theorem may be generalized to the case when $M$ has more
serious singularities (say, to the case when $M$ is a singular algebraic
manifold). In this case, the usual cohomology must be replaced by the
intersection cohomology. The use of the intersection cohomology also allows to
release the assumption that 0 is a regular value of the moment map. 

d. \ The special case of the theorem, when $M$ is a toric manifold, is due to
R.~Stanley \cite{Stanley87}. In this sense, our result is a generalization of
Stanley's theorem. In particular, we obtain a new, more geometric proof, of the
Stanley's theorem. See \refss{Stanley} for details.
\erem

\ssec{sketch}{}
The proof of the theorem is based on a study of the action of $\tet$ on the
equivariant cohomology of $M$. We now formulate the main results about this
action. The proofs are given on Sections~\ref{S:HT=HxA} and \ref{S:Morse}. Some
examples and applications of \reft{main} are discussed in \refs{Ex}.

\ssec{eqcoh}{Action of the involution on the equivariant cohomology}
Let $ET$ denote the universal $T$-bundle and let $BT=ET/T$ denote the
classifying space of $T$. Let $H_T^*(M)=H^*(ET\times_TM;\CC)$ denote the
equivariant cohomology of $M$ with complex coefficients. Then $\tet$ acts
naturally on $H_T^*(M)$, cf. \refss{action}. Let
\eq{chi}
     \chi_\tet(t) \ := \ \sum_{t=0}^\infty\, t^i\Tr \tet|_{H_T^i(M)}
\end{equation}
denote the {\em graded character} of this action.

In \refs{HT=HxA} we use the fact that the spectral sequence of the fibration 
$M\times_TET\to BT$ degenerates at the second term (cf. \cite[Proof of
Pr.~5.8]{Kirwan84}, \cite[Proof of Th.~5.3]{BrFar3}) to prove the following
proposition, in which we do not assume that 0 is a regular value of the moment
map.

\prop{HT=HxA}
\(
\displaystyle
        \chi_\tet(t) \ = \
                \sum_{i=1}^{2n} (h^{i,+} - h^{i,-}) \frac{t^{i}}{(1+t^2)^k}.
\)
\eprop

On the other side, in \refs{Morse} we use a version of the equivariant Morse
inequalities \cite{AtBott82} constructed in \cite{BrFar4,BrFar3} for functions non-degenerate in the sense of Bott and extended to Kirwan-Morse functions  in \cite{BrSil}  to get the
following

\prop{chitet}
Suppose that zero is a regular value for the moment map $\mu$ and let
$M_0:=\mu^{-1}(0)/T$. Then $\chi_\tet$ is equal to the graded character of the
action of $\tet$ on $H^*(M_0)$:
\eq{chitet}
   \chi_\tet(t) \ =  \ \sum_{i}\,  (h^{i,+}_0 - h^{i,-}_0) t^i.
\end{equation}
\eprop
\cor{toric}
If, in the conditions of \refp{chitet}, the dimension of the torus is equal to
$n=\frac12\dim_\RR{M}$, then $\chi_\tet=1$.
\ecor

\ssec{prmain}{Proof of the main theorem}
Comparing Propositions~\ref{P:HT=HxA} and \ref{P:chitet} we obtain
\[
        \sum_{i=1}^{2n} (h^{i,+} - h^{i,-}) \frac{t^{i}}{(1+t^2)^k} \ = \
                \sum_{i}\, (h^{i,+}_0 - h^{i,-}_0) t^i,
\]
which is equivalent to \reft{main}.
\hfill$\square$


\sec{Ex}{Examples and applications}

\ssec{Stanley}{Symmetric toric variety. Application to combinatorics}
\reft{main} takes a particularly simple form if the dimension of the torus is equal to
$n=\frac12\dim_\RR M$, so that $M$ is a toric variety. In this case we
say that $M$ is a {\em symmetric (with respect to the involution
$\tet$) toric variety}.

If $M$ is a symmetric toric variety, then
the reduced space $M_0$ is a point. Hence, $h^{i,-}_0=0$ for all $i$,
$h^{i,+}_0=0$ for all $i>0$ and $h^{0,+}_0=1$. Thus \reft{main} reduces
in this case to the following statement, which was originally proven
by R.~Stanley by a completely different method
\footnote{Stanley proved the theorem for more general case, when $M$
  is a symmetric toric orbifold. One can prove this result using our
  method and \refr{main}.b}.
\cor{Stanley} If, in the conditions of \reft{main}, the dimension of the
torus is equal to $n=\frac12\dim_\RR M$, then
\[
        h^{i,+} \ - \ h^{i,-} \ = \ \binom{n}{i}.
\]
\ecor
There are a lot of examples of symmetric toric varieties. To describe these
examples let us recall that each toric variety is completely characterized by
its {\em moment polytope} $\mu(M)\subset \grt^*$. The toric variety $M$ is an
orbifold if and only if the polytope $\mu(M)$ is {\em simple}, i.e. if each of
its vertices has the valence $n=\frac12\dim_\RR M$. There is a complete
description of polytopes corresponding to smooth toric varieties, cf., for
example, \cite[\S{}IV.2]{Audin91}. The equation \refe{normal} implies that
the moment polytope corresponding to a symmetric toric variety is centrally
symmetric. Vice versa, if the moment polytope is centrally symmetric, one easily
constructs an involution on $M$ satisfying \refe{normal}. Hence, symmetric toric
orbifolds are in one-to-one correspondence with centrally symmetric simple
convex polytopes. \refc{Stanley} can be used to get an estimate on the number of
faces of such a polytope. See \cite{Stanley87} for details.
\rem{Campo} A. A\'{}Campo-Neuen \cite{Campo} extended
\refc{Stanley} to singular toric varieties. This leads to an extension of
Stanley's estimates on the number of faces of a centrally symmetric polytopes to
rational polytopes, which are not necessarily simple. This result may be also
obtained by our method, cf. \refr{main}.c.
\erem

\ssec{CPn}{$\CC{P}^3$ as an $S^1$-space with involution}
Some toric varieties, which are not symmetric, still posses an
involution compatible, in the sense of \refe{compat}, with the action
of a torus of smaller dimension. Before discussing more general
examples, let us consider the simplest case $M=\CC{P}^{3}$. Then the
moment polytope is a $3$-dimensional simplex, which is, obviously, not
centrally symmetric. Hence, there is no involution on $M$ compatible,
in the sense of \refe{compat}, with the action of a $3$-dimensional
torus. However, using the homogeneous coordinates $[z_1:z_2:z_3:z_4]$
on $M$, we can  define the action of the circle $T=S^1$ on $M$ and an
involution $\tet:M\to M$ by
\begin{gather}
    t \ \cdot \ [z_1:z_2:z_3:z_4] \ = \ [tz_1:tz_2:z_3:z_4], \qquad
                        t\in S^1=\{\zet\in\CC:\, |\zet|=1\} \notag \\
    \tet \ \cdot \ [z_1:z_2:z_3:z_4] \ = \ [z_3:z_4:z_1:z_2].\notag
\end{gather}
Clearly, all the conditions of \reft{main} are satisfied.

In this case, both sides of \refe{main} may be easily calculated. In particular,
$M_0=\CC{P}^1\times\CC{P}^1$ and $\tet$ acts on $M_0$ by
$\tet:(a,b)\mapsto{}(b,a)$, where $a,b\in\CC{P}^1$. An interested
reader can easyly verify \reft{main} in this simple case.

\ssec{sym}{Toric varieties with involution}
The previous example allows the following generalization. Let $M$ be a
toric variety endowed with an action of the torus $T$ of dimension
$n=\frac12\dim_\RR M$. Suppose that there exists a linear involution
$\tiltet:\grt^*\to\grt^*$ which preserves the moment polytope
$\mu(M)$.  Then $\tiltet$ induces an involution $\tet:M\to M$.
However,  in general, this involution is not compatible with the
action of $T$ in the sense of \refe{compat}.

Let $\tiltet^*:\grt\to\grt$ denote the involution of the Lie algebra $\grt$ dual
to $\tiltet$. Let $T'\subset T$ be a subtorus, such that $\tiltet^*$ acts as
multiplication by $-1$ on the Lie algebra of $T'$. Then the actions of $T'$ and
$\tet$ on $M$ are compatible in the sense of \refe{compat}. If, in addition, 0
is a regular value of the moment map for the $T'$ action, then \reft{main} may
be applied.

\ssec{flag}{The variety of complete flags in $\CC^n$}
We now give an example of a manifold, which is not toric, but
satisfies the conditions of \reft{main}.

Let $\lam=\{\lam_1\nek\lam_n\}$ be a centrally symmetric set of real numbers. In
other words, $\lam=\{\pm\nu_1\nek\pm\nu_k\}$ if $n=2k$ and
$\lam=\{0,\pm\nu_1\nek\pm\nu_k\}$ if $n=2k+1$. We will assume that the numbers
$\nu_i$ above are positive and mutually different. The set $A_\lam$ of all complex
Hermitian $n\times{}n$-matrices with spectrum $\lam$ is naturally identified
with the variety of complete flags in $\CC^n$. In particular, $A_\lam$ has a
structure of a compact K\"ahler manifold.

Fix mutually different rational numbers $r_1\nek r_n$ and consider the action of
the circle $S^1=\{e^{it}:\, t\in\RR\}$ on $A_\lam$, defined by
\eq{actionF}
    e^{it}\cdot A \ = \
       \diag\big\{\, e^{ir_1t}\nek e^{ir_nt}\, \big\} \ \cdot \ A \
       \cdot \
          \diag\big\{\, e^{-ir_1t}\nek e^{-ir_nt}\, \big\},
       \qquad t\in\RR, \ A\in A_\lam.
\end{equation}
This action is Hamiltonian with respect to the \ka
structure on $A_\lam$: if we identify the coalgebra Lie of $S^1$ with
$\RR$, then the function
\[
      \mu:\, A_\lam \ \to \RR, \qquad A=\{a_{ij}\} \ \mapsto \
      \sum_{i=1}^n r_ia_{ii}
\]
is a moment map for this action.

Define an involution $\tet:A_\lam\to A_\lam$ by the formula
\[
       \tet:\, A \  \mapsto \ -A^t.
\]
This involution preserves the \ka structure on $A_\lam$ and satisfies
\refe{compat}, \refe{normal}. Moreover, one easily checks that $S^1$ acts
locally freely on $\mu^{-1}(0)$. Hence, zero is a regular value of $\mu$ and
\reft{main} is applicable.

\rem{flag}
a. \ The action \refe{actionF} is a restriction of the conjugate action on
$A_\lam$ of the torus of all unitary diagonal matrices. The latter action is
also Hamiltonian and compatible with the involution $\tet$ in the sense of
\refe{compat}. Unfortunately, zero is not a regular value of the moment map for
this action. However, the generalization of \reft{main}, indicated in
\refr{main}.c, may be applied in this case.

b. \ One easily generalizes the results of this subsection to a flag variety of
an arbitrary reductive group $G$. The involution $\tet$ should be replaced by
the action of the longest element of the Weyl group. Also the action of the
maximal torus of $G$ should be restricted to a subtorus on whose Lie algebra the
longest element of the Weyl group acts as a multiplication by $-1$, cf.
\refss{sym}. We leave the details to an interested reader.
\erem

\sec{HT=HxA}{Proof of \refp{HT=HxA}}

\ssec{action}{Action of $\tet$ on $H_T^*(M)$}
Before proving the proposition let us describe more explicitly the
action  of $\tet$ on $H_T^*(M)$. One of the ways to
define this action is the following. Consider the group
\[
        \tilT \ := \ \big\{\, (\eps,t):\, \eps=\pm1, t\in T\, \big\}
\]
with the product law
\[
        (\eps_1,t_1)\cdot (\eps_2,t_2) \ = \ (\eps_1\eps_2, t_1^{\eps_2}t_2).
\]
Then $\tilT$ acts on $M$ by the formula
\[
        (\eps,t)\cdot m  \ = \ t^\eps\cdot \tet^{\frac{1-\eps}2}(m).
\]
We will identify $\tet$ with the element $(-1,1)\in \tilT$. Let
$E\tilT$ denote the universal $\tilT$-bundle. Then $\tet$ acts
diagonally on $E\tilT\times_TM$ and, hence, on
$H_T^*(M)=H^*(E\tilT\times_TM)$ (in the last equality we used that
$E\tilT$ may be considered also as a model for the universal
$T$-bundle).

\ssec{tetBT}{Action of the involution on $BT$}
We consider the quotient $BT={E\tilT}/T$ as a model for the classifying space of
$T$. Then $\tet=(-1,1)\in\tilT$ acts on $BT$ and, hence, on the cohomology ring
$H^*(BT)$. The later is isomorphic to the graded ring $\CC[\grt]$ of polynomials
on the Lie algebra $\grt$ of $T$ (with grading given by twice the degree of the
polynomial):
\eq{isom}
       H^*(BT) \ \simeq\CC[\grt].
\end{equation}
The following lemma describes the action of $\tet$ on this ring.
\lem{tetBT}
Under the isomorphism \refe{isom} the action of $\tet$ on $H^*(BT)$ is
given by the rule
\[
        \tet(t) \ = \ -t, \qquad t\in \grt.
\]
In particular, the restriction of $\tet$ to $H^2(BT)$ is equal to $-1$.
\elem
\prf
This lemma is well known, but in view of the difficulty we have to
locate an explicit reference, we present a proof here.
For simplicity, we only present the proof for the case $\dim_\RR{}T=1$.
The arguments in general case are exactly the same, but the notation
is more complicated. We also identify the one dimensional torus with
the unit circle $S^1=\{e^{i\phi}:\, \phi\in\RR\}$ in $\CC$.

Fix $m\ge3$ and consider the free action of $\tilT$ on the product of spheres
\[
      S^{2m-1}\times S^{2m-1} \ := \ \big\{\,
              (z_1\nek z_m;w_1\nek w_m)\in \CC^{2m}:\,
                |z_1|^2+\cdots|z_m|^2=1;\, |w_1|^2+\cdots|w_m|^2=1
                       \big\},
\]
given by
\begin{gather}
      e^{i\phi}\cdot(z_1\nek z_m;w_1\nek w_m) \ \mapsto
        (e^{i\phi}z_1\nek e^{i\phi}z_m;w_1\nek w_m); \notag
      \\
      (-1,0)\cdot(z_1\nek z_m;w_1\nek w_m) \ \mapsto
        (\oz_1\nek \oz_m;-w_1\nek -w_m). \notag
\end{gather}
Let us denote by $BT_m$ the quotient of $S^{2m-1}\times S^{2m-1}$ by the action
of $T\subset\tilT$. Since the reduced cohomology of $S^{2m-1}\times S^{2m-1}$
vanishes up to dimension $2m-2\ge4$ it follows that there is a natural $\tilT/T$
equivariant isomorphism between the cohomology $H^i(BT_m)$ and $H^i(BT)$ for any
$i$ less than $2m-3\ge 3$ (this may be shown by the same arguments as in the
proof of Lemma~2.8 in \cite{BrFar3})

Clearly, $\tet$ acts on $BT_m=\CC\PP^{m-1}\times S^{2m-1}$ as complex
conjugation on the first factor and multiplication by $-1$ on the second factor.
Hence, the induced action on the cohomology $H^2(BT_m)$ is multiplication by
$-1$. This proves the lemma.
\eprf
\cor{charBT}
The graded character of the action of $\tet$ on $H^*(BT)$ is equal to
$(1+t^2)^{-k}$.
\ecor

\ssec{spseq}{The spectral sequence. Proof of \refp{HT=HxA}}
The projection $E\tilT\times M\to E\tilT$ induces a $\tilT/T$
equivariant fiber bundle $p:E\tilT\times_TM\to BT$. The fiber of this
bundle is $\tilT/T$ equivariantly homeomorphic to $M$.  The cohomology
$H^*(E\tilT\times_TM)=H^*_T(M)$ may be calculated by the spectral
sequence of the above bundle. The later spectral sequence degenerates
at the second term  (cf. \cite[Proof of
Pr.~5.8]{Kirwan84}, \cite[Proof of Th.~5.3]{BrFar3}) and is obviously
$\tilT/T$ invariant. It follows that there exists a $\tilT/T$
equivariant isomorphism of graded rings
\[
       H^*_T(M) \ \simeq \ H^*(M) \otimes_\CC H^*(BT).
\]
Hence, the graded character $\chi_\tet(t)$ of the action of $\tet$ on $H^*_T(M)$
is equal to the product of the graded characters of the actions of $\tet$ on
$H^*(M)$ and $H^*(BT)$. \refp{HT=HxA} follows now from \refc{charBT}.
\hfill$\square$


\sec{Morse}{Equivariant Kirwan-Morse inequalities and proof of \refp{chitet}}

Our proof of \refp{chitet} is an application of the
$\tilT$-equivariant Morse inequalities for the square of the moment
map. Since the group $\tilT$ is disconnected, it is important to use
the equivariant Morse inequalities with local coefficients  (cf. \cite{BrFar4,BrFar3} for the Morse-Bott case and  \cite[Th.~2.6]{BrSil} for the general case).

The key fact of the proof is that the square of the moment map is an
{\em equivariantly perfect Morse function}, i.e. the above inequalities
are, in fact, equalities, cf. \refl{perfect}. This result is a slight
generalization of a theorem of Kirwan \cite[Ch~6, Th.~5.4]{Kirwan84}, who
considered $T$-equivariant Morse inequalities with trivial local
coefficients.

\ssec{morse}{Equivariant Morse inequalities}
For convenience of the reader we recall here equivariant Morse inequalities with
local coefficients for an action of a disconnected group, as they formulated in
\cite{BrFar4,BrFar3,BrSil}. In this subsection $M$ is a compact manifold acted upon by
a compact not necessarily connected Lie group $G$. Let $\calF$ be a
$G$-equivariant flat vector bundle over $M$. Denote by $H^*_G(M,\calF)$ the
$G$-equivariant cohomology of $M$ with coefficients in $\calF$ and let
\[
    \calP^G_\calF(t) \ = \ \sum_{i}\, t^i \dim_\CC\, H^i_G(M,\calF)
\]
be the {\em equivariant Poincar\'e series} of $M$ with coefficients on $\calF$.

Recall that a smooth function $f:M\to\RR$ is called  {\em Kirwan-Morse} if it satisfies the following non-degeneracy conditions:
\begin{itemize}
\item[\textbf{(C1)}]
The critical set $C$ is a finite union of disjoint closed subsets $Z\subset  C$ called {\em
components of the critical set of $f$} such that $f|_Z=\const$.

\item[\textbf{(C2)}]  For component $Z\subset C$ there exists a locally closed connected submanifold
$\Sig_Z$, called a {\em minimizing manifold}, containing $Z$ such that for every $x\in \Sig_Z, \ x\not\in{}Z$ we have 
$f(x)>f|_Z$.  In other words, $Z$ is the subset of $\Sig_Z$ on  which $f$ takes its minimum values.

\item[\textbf{(C3)}]
At every point $x\in C$ the tangent space $T_x\Sig_Z$ is the maximal among all subspaces of
$T_xM$ on which the Hessian $H_x(f)$ is positive semi-definite.
\end{itemize}
Thus the  critical subsets of a Kirwan-Morse function can be as degenerate as a minimum of a
function, but not worse.

Suppose $f:M\to\RR$ is a $G$-invariant Kirwan-Morse function on $M$. Let $Z$ be a component of the critical point set $C$ and let
$\nu^-(Z)$ denote the normal bundle to $\Sig_Z$ in $M$.   The dimension of $\nu^-(Z)$ is called the {\it index} of $Z$ (as a component of the critical set of $f$) and is denoted by $\ind(Z)$.  Let $o(Z)$ denote the {\it orientation bundle of $\nu^-(Z)$, considered as a flat line
bundle}.

If the group $G$ is connected, then $Z$ is a $G$-invariant subset of $M$. In general, we denote by 
\[
	G_Z\ =\ \big\{g\in G| \ g\cdot Z\subset Z\big\}
\] 
the stabilizer of the component $Z$ in $G$.  Let $|G:G_Z|$ denote the index of $G_Z$ as a subgroup of $G$; it is always finite.

The compact Lie group $G_Z$ acts on $Z$ and the flat vector bundles $\F{|_Z}$ and $o(Z)$ are $G_Z$-equivariant.  Let
 \(
        \check{H}_{G_Z}^\ast(Z,\fo)
 \)
denote the {\em equivariant \Cech cohomology} of the flat $G_Z$-equivari\-ant vector bundle
$\fo$.  Consider the {\it
equivariant Poincar\'e series}
$$
        \calP^{G_Z}_{Z,\F}(t)\ =\
           \sum_i t^i \dim_{\CC}  \check{H}_{G_Z}^i(Z,\fo)
$$
and define, using it, the following {\it equivariant Morse counting
series}
$$
        \calM_{f,\F}^G(t)\ = \ \sum_Z
             t^{\ind(Z)} |G:G_Z|^{-1}\calP_{Z,\F}^{G_Z}(t)
$$
where the sum is taken over all components $Z$ of $C$.

The following version of the equivariant Morse inequalities is a
particular case of \cite[Th.~2.6]{BrSil} (see also  \cite[Th.~7]{BrFar4},\cite[Th.~1.7]{BrFar3}).
\th{BrFar}
   Suppose that $G$ is a compact  Lie group, acting on a closed manifold
   $M$, and let $\F$ be an equivariant flat vector bundle over $M$.
   Then for any   $G$-equivariant Kirwan-Morse
   function
   $f:M\to\RR$, there exists a formal power series $\calQ(t)$ with
   non-negative integer coefficients, such that
   \[
        \calM_{f,\F}^G(t) \ - \ \calP_{\F}^G(t) \ = \ (1+t)\calQ(t).
   \]
\eth

\ssec{eqbundle}{Equivariant flat bundles}
We now return to the situation described in \refs{main}.
Let $\calF$ be a $\tilT$-equivariant flat vector bundle over $M$. This
is the same as $T$-equivariant flat vector bundle, on which $\tet$ acts
preserving the connection and so that $\tet(t\cdot \xi)=(t^{-1})\cdot\tet(\xi)$ for
any $\xi\in\calF, t\in T$.

In our proof of \refp{chitet} we will only use the following type of
$\tilT$-equivariant bundles:  Let $\rho:\ZZ_2\to \End_\CC(V_\rho)$ be
a representation of $\ZZ_2$ in a finite-dimensional complex vector
space $V_\rho$. Consider the bundle $\calF_\rho=M\times V_\rho$ with
the trivial connection and with the action of $\tilT$ given by
\[
     (\eps,t):\, (m,\xi) \mapsto \big(\, (\eps,t)\cdot m,\, \rho(\eps)\cdot\xi\,
     \big).
\]
Then the the equivariant cohomology $H_{\tilT}^*(M,\calF_\rho)$ of $M$ with
coefficients in $\calF_\rho$ is given by
\eq{HFrho}
   H_{\tilT}^*(M,\calF_\rho) \ = \
       \Hom_{\ZZ_2}\big(V_\rho^*, H_{T}^*(M)\big)
\end{equation}
(here $V^*_\rho$ is the representation dual to $V_\rho$). In particular, if
$\rho$ is the regular representation of $\ZZ_2$ then
\eq{reg}
        H_{\tilT}^*(M,\calF_\rho)= H_{T}^*(M).
\end{equation}

The bundle $\calF_\rho$ is completely determined by the representation $\rho$.
When it causes no confusion, we will write $\rho$ for $\calF_\rho$ in order to
simplify the notation.

\ssec{f}{Morse equalities for the square of the moment map}
Fix a $T$-invariant scalar product $\<\cdot,\cdot\>$ on $\grt^*$ and consider the real-valued
function $f=\<\mu,\mu\>$ on $M$. By \cite[Ch.~I,  \S4]{Kirwan84} this is a Kirwan-Morse function. Since $\mu$ is $T$-invariant so is $f$. Also  $\mu\circ\tet=- \mu$ implies that $f\circ\tet= f$. We conclude that $f$ is a $\tilT$ invariant Kirwan-Morse function. Moreover, cf, \cite[Ch.~I,\S4.18]{Kirwan84},  if $Z$ is  a component of the set $C$ of critical
points of $f$ (cf. \refss{morse}), then $\nu^{-}(Z)$  is even dimensional $T$-equivariant bundle. The action of $T$ defines an orientation of $\nu^-(Z)$ and, hence, the orientation bundle $o(Z)$ is trivial.

The set $Z_0:=\mu^{-1}(0)$ is a connected component of the set of critical
points of $f$. The group $\tilT$ acts on $Z_0$ and the equivariant
cohomology
\[
       H^*_\tilT(Z_0,\calF_\rho|_{Z_0}) \ = \    \Hom_{\ZZ_2}\big(V^*_\rho,H^*_T(\mu^{-1}(0))\big).
\]
Recall that $T$ acts locally freely on $\mu^{-1}(0)$ and we denote $M_0=\mu^{-1}(0)/T$. Hence, $H^*_T(\mu^{-1}(0))= H^*(M_0)$ and we get
\[
       H^*_\tilT(Z_0,\calF_\rho|_{Z_0}) \ = \  \Hom_{\ZZ_2}\big(V^*_\rho,H^*(M_0)\big).
\]
The equivariant Poincar\'e series of $Z_0$ with coefficients in
$\calF_\rho$ is given by
\[
        \calP^\tilT_{Z_0,\rho}(t) \ = \
            \sum_i\, t^i\dim_\CC H^*_\tilT(Z_0,\calF_\rho|_{Z_0})
               \ = \
                \sum_i\, t^i\dim_\CC\Hom_{\ZZ_2}\big(V^*_\rho,H^*(M_0)\big).
\]
In particular,
\eq{Z0}\begin{aligned}
   \calP^\tilT_{Z_0,\rho}(t) \ &= \ \sum_i\, t^ih^{i,+}_0,
   \quad\text{if $\rho$ is the trivial representation}; \\
   \calP^\tilT_{Z_0,\rho}(t) \ &= \ \sum_i\, t^ih^{i,-}_0,
   \quad\text{if $\rho$ is the sign representation}.
\end{aligned}\end{equation}
The stabilizer $G_{Z_0}$ of $Z_0$ coincides with the whole group $\tilT$. Thus
$|G:G_{Z_0}|=1$.

The involution $\tet$ acts freely on the set of  components of $C$
different from $Z_0$, cf. \refr{main}.a. Hence, if $Z\not=Z_0$ is a component of $C$, then $G_Z=T$ and $|G:G_Z|=2$. The $T$-equivariant Poincar\'e
polynomial of $Z$ depends only on the dimension of $\rho$ and is given by
\[
      \calP^T_{Z}(t)   \ = \ \dim_\CC V_\rho\cdot   \sum_i\, t^i\dim_\CC \check{H}^i_T(Z).
\]
Hence, the equivariant Morse counting series
\[
       \calM^{\tilT}_{f,\rho}(t) \ = \ \calP^\tilT_{Z_0,\rho}(t) \ + \ \frac12\sum_Z\, \calP^T_{Z}(t).
\]
Here the sum in the right hand side is taken over all components $Z$ of the set of critical points of $f$ different from $Z_0$.

Let
\[
       \calP^\tilT_\rho(t) \ = \ \sum_i\, t^i \dim_\CC H^i_\tilT(M,\calF_\rho)
\]
be the equivariant Poincar\'e series of $M$ with coefficients in $\calF_\rho$.

The following lemma expresses the fact that $f$ is an {\em equivariantly perfect Morse function}.
\lem{perfect}
The following equality holds
\eq{perfect}
        \calM^{\tilT}_{f,\rho}(t) \ = \ \calP^\tilT_\rho(t).
\end{equation}
\elem
\prf
It follows from \reft{BrFar}, that there exists a formal power series
$Q_\rho(t)$ with non-negative coefficients, such that
\eq{morse}
        \calM^{\tilT}_{f,\rho}(t) \ = \ \sum_i\, t^i H^i_\tilT(M,\calF_\rho) \ + \ (1+t)\,Q_\rho(t).
\end{equation}
Our goal is to show that $Q_\rho\equiv0$.

It follows from \refe{HFrho}, that both the equivariant Morse
counting series and the equivariant Poincar\'e series are additive with respect
to $\rho$. More precisely, if $\rho_1\oplus\rho_2$ denotes the direct sum of two
representations then
\[
        \calM^{\tilT}_{f,\rho_1\oplus\rho_2}(t) \ = \
                \calM^{\tilT}_{f,\rho_1}(t) \ + \ \calM^{f,\tilT}_{\rho_2}(t), \quad
       \calP^\tilT_{\rho_1\oplus\rho_2}(t) \ = \
                \calP^\tilT_{\rho_1}(t) \ + \ \calP^\tilT_{\rho_2}(t).
\]
Hence, it suffices to prove the lemma for the irreducible representations of
$\ZZ_2$. Moreover, it follows from \refe{morse} that it is enough to prove that
$Q_\rho=0$ when $\rho$ is a reducible representation, which contains any of the
irreducible representation as a subrepresentation.%

In particular, it is enough to prove that $Q_\rho=0$ when $\rho$ is
the regular representation. However, \refe{reg} implies that, if
$\rho$ is the regular representation, then \refe{morse} reduces to the
$T$-equivariant Kirwan-Morse inequalities with trivial coefficients. It was
shown by Kirwan \cite[Th.~5.4]{Kirwan84} that the later inequalities
are exact, i.e. $Q_\rho=0$.  \eprf

\ssec{proof}{Proof of \refp{chitet}}
Set
\[
        H^*_T(M)^\pm \ := \ \big\{\, x\in H^*_T(M):\, \tet x=\pm x\, \big\}
\]

Let us, first, apply \refl{perfect} with $\rho$ being the trivial representation
of $\ZZ_2$. It follows from \refe{HFrho} that, in this case,
$\calP^\tilT_{\rho}(t)=\sum_i\, t^i \dim_\CC H^i_T(M)^+$. Hence, from
\refe{Z0} and \refl{perfect} we obtain
\eq{1}
         \sum_i\, t^i \dim_\CC H^i_T(M)^+ \ = \
           \frac12\sum_Z\, \calP^T_{Z}(t) \ + \
             \sum_i\, t^ih^{i,+}_0.
\end{equation}

Apply now \refl{perfect} with $\rho$ being the sign representation of $\ZZ_2$.
Then $\calP^\tilT_{\rho}(t)=\sum_i\, t^i \dim_\CC H^i_T(M)^-$ and
$\calP^\tilT_{Z_0,\rho}=\sum_i\, t^ih^{i,-}_0$. Hence,
\eq{2}
        \sum_i\, t^i \dim_\CC H^i_T(M)^- \ = \
                \frac12\sum_Z\, \calP^T_{Z}(t) \ + \
             \sum_i\, t^ih^{i,-}_0.
\end{equation}
Subtracting \refe{2} from \refe{1} we get \refe{chitet}.
\hfill$\square$

\providecommand{\bysame}{\leavevmode\hbox to3em{\hrulefill}\thinspace}

\end{document}